\documentclass[11pt]{amsart}%
\usepackage{amssymb}
\usepackage{amsmath}
\usepackage{amsfonts}
\usepackage{graphicx}
\usepackage{tabulary}%
\newtheorem{theorem}{Theorem}[section]
\newtheorem{corollary}[theorem]{Corollary}
\newtheorem{lemma}[theorem]{Lemma}
\newtheorem{proposition}[theorem]{Proposition}

\newtheorem{remark}[theorem]{Remark}

\begin{document}
\title{Uniqueness of static decompositions}
\author[M. Guti\'errez]{Manuel Guti\'errez}
\address{M. Guti\'errez. Departamento de Algebra, Geometr\'ia y Topolog\'ia.
Universidad de M\'alaga. Spain.}
\email{mgl@agt.cie.uma.es}
\author[B. Olea]{Benjam\'in Olea}
\address{B. Olea. Departamento de Algebra, Geometr\'ia y Topolog\'ia. Universidad de
M\'alaga. Spain.}
\email{benji@agt.cie.uma.es}
\thanks{The first author was supported in part by MEYC-FEDER Grant MTM2007-60016.}

\begin{abstract}
We classify static manifolds which admit more than one static decomposition
whenever a condition on the curvature is fullfilled. For this, we take a
standard static vector field and analyze its associated one parameter family
of projections onto the base. We show that the base itself is a static
manifold and the warping function satisfies severe restrictions, leading us to
our classification results. Moreover, we show that certain condition on the
lightlike sectional curvature ensures the uniqueness of static decomposition
for Lorentzian manifolds.

\end{abstract}
\maketitle



\section{Introduction}

Given $(L,g_{L})$ a connected Riemannian manifold, $f\in C^{\infty}(L)$ a
positive function and $\varepsilon=\pm1$, we call \textit{static manifold} to
a product $L\times\mathbb{R}$ furnished with the metric $g_{L}+\varepsilon
f^{2} dt^{2}$, which is denoted by $L\times_{\varepsilon f}\mathbb{R}$. To
follow standard definitions \cite{ONeill}, we refer to the case $\varepsilon
=-1$ as \textit{static space} instead of static manifold. Although static
manifold is not a common term, it allows us to handle jointly the Riemannian
and the Lorentzian case.

A vector field is irrotational if it has integrable orthogonal distribution.
It is well known that a manifold furnished with an irrotational and Killing
vector field can be decomposed as a static manifold, at least locally. Due to
this, such a vector field is called \textit{static} and when it gives rise to
a global decomposition, it is called \textit{standard static}, \cite{Sanchez1}. Obviously, the existence of two different static decompositions is
equivalent to the existence of two standard static vector fields linearly
independent at some point. Minkowski, euclidean, hyperbolic and a suitable portion of the anti De
Sitter spaces provide us good examples of manifolds with different decompositions, but there are many other manifolds which can be also decomposed as a static manifold in several ways.

The decomposition uniqueness problem has been studied by several authors. The
De Rham-Wu Theorem ensures the uniqueness of a direct product decomposition
for a simply connected semi-Riemannian manifold (for the nonsimply connected
case see \cite{Esche} and \cite{Gut2}) and the uniqueness of Generalized
Robertson-Walker decomposition was studied in \cite{Gut}, obtaining that the
De Sitter space is the only complete space with several nontrivial
decomposition, whereas Friedmann spaces have a unique decomposition even locally.

The uniqueness of static decomposition seems more complicated and only partial
results have been obtained. In \cite{Sanchez1} it is shown that static spaces
with compact base do not admit another global static decomposition and in
\cite{Tod} it is computed Riemannian three dimesional metrics which can be
used to construct Einstein static spacetimes in more than one way. On the
other hand, it is a remarkable fact that for the exterior Schwarzschild
spacetime the uniqueness can be directly shown. Indeed, any timelike Killing
vector field is proportional to the canonical one, \cite{ONeill}. Apart from
these results, no much more is known about this topic.

In this paper we study manifolds with more than one static decomposition, with
special emphasis on the Lorentzian case. For this, we consider a static
manifold with a standard static vector field linearly independent to the
canonical one at some point. The projection onto the base gives us a family of
vector fields, which we use to decompose the base itself as a static manifold.
This allows us to prove in Proposition \ref{teoclasifica} that these manifolds
are a special type of warped product and, under a mild curvature hypothesis
about the base or Einstein assumption, we classify them in Theorem
\ref{ClasificaCurvatura} and \ref{Einstein} respectively. We particularize to
Einstein spacetimes in Theorem \ref{Einsteinspacetime} and finally, in Theorem
\ref{Unicidad}, we show that in the Lorentzian case the uniqueness is
guaranteed if the lightlike sectional curvature at a point is never zero.

\section{Preliminaries}

Given a product manifold $L\times\mathbb{R}$, the lift of $\partial_{t}$ is
still denoted by $\partial_{t}$, however, given $X\in\mathfrak{X}(L)$ we will
denote by $\tilde{X}$ to its lift. We call $\pi:L\times\mathbb{R}\rightarrow
L$ the canonical projection and $i_{t}:L\rightarrow L\times\mathbb{R}$ the
injection given by $i_{t}(p)=(p,t)$. We will usually avoid $\pi$ in the
formulaes to lighten the notation.

If $V\in\mathfrak{X}(L\times\mathbb{R})$ we call $V_{p}^{t}=\pi_{\ast
(p,t)}(V)$, which is a vector field on $L$ for each $t\in\mathbb{R}$. Fixed
$p\in L$, $V^{t}_{p}$ is a curve in $T_{p}L$ and thus we can consider the
vectorial derivative $\frac{d}{dt}V_{p}^{t}\in T_{p}L$, obtaining in this way
another vector field on $L$ for each $t\in\mathbb{R}$. On the other hand, if
$h\in C^{\infty}(L\times\mathbb{R})$ then we call $h^{t}$ the function on $L$
given by $h^{t}=h\circ i_{t}$, but $h_{t}$ will mean the derivative respect to
$t$.

Recall that any vector field on a two or one dimensional manifold is
irrotational. By convention, in a one dimensional manifold, the orthogonal
leaf of a vector field through a point simply means this point. This will allow
us to prove some results without distinguishing cases.

We denote by $\mathbb{R}_{[\varepsilon]}^{2}$ the Euclidean/Minkowski plane
$(\mathbb{R}^{2},ds^{2}+\varepsilon dt^{2})$ and $\mathbb{H}_{[\varepsilon
]}^{2}(r)$ the hyperbolic/anti De Sitter plane $(\mathbb{R}^{2},ds^{2}%
+\varepsilon\cosh^{2}(rs)dt^{2})$ of curvature $-r^{2}$. By $\widehat
{\mathbb{H}}_{[\varepsilon]}^{2}(r)$ we denote the static manifold
$(\mathbb{R}^{2},ds^{2}+\varepsilon e^{2rs}dt^{2})$, which is another
representation of the hyperbolic plane if $\varepsilon=1$ and a piece of the
anti De Sitter plane if $\varepsilon=-1$. To simplify this notation, we also
denote $\mathbb{H}_{[1]}^{2}(r)$ by $\mathbb{H}^{2}(r)$ and $\mathbb{R}%
_{[1]}^{2}$ by $\mathbb{R}^{2}$.

If $\Pi$ is a degenerate plane in a Lorentzian manifold of dimension greater
than two, the lightlike sectional curvature of $\Pi$ is defined as
\[
\mathcal{K}_{u}(\Pi)=\frac{g(R(v,u,u),v)}{g(v,v)},
\]
where $u,v\in\Pi$ with $u$ lightlike and $v$ spacelike. This curvature depends
on the choosen $u$, but its sign does not depend on it. Thus, we can say zero
lightlike sectional curvature without explicit mention of the choosen
lightlike vector. We can easily compute the lightlike sectional curvature in a
static space.

\begin{lemma}
\label{curluz} Let $M=L\times_{-f}\mathbb{R}$ be a static space with $dim
L\geq2$ and take $v,w\in TL$ unitary vectors with $v\perp w$. If
$\Pi=span(v,u)$, where $u=w+\frac{1}{f}\partial_{t}$, then
\[
\mathcal{K}_{u}(\Pi)=K^{L}\left(  span(v,w)\right)  +\frac{g(\nabla
_{v}\widetilde{\nabla f},v)}{f},
\]
being $K^{L}$ the sectional curvature of $L$.
\end{lemma}

The following result is well known. We prove the precise version that we need which,
jointly with the below remark, ensures the global decomposition as a static
manifold.

\begin{proposition}
\label{covering} Let $M$ be a geodesically complete Lorentzian or Riemannian manifold of dimension
greater that one. If $V$ is a static vector field without zeros (timelike in
the Lorentzian case), then there exists a normal semi-Riemannian covering map
$\Phi:\left(  L\times\mathbb{R},g_{L}+\varepsilon f^{2}dt^{2}\right)
\rightarrow(M,g)$ with $\Phi_{*}(\partial_{t})=V$.
\end{proposition}

\begin{proof}
Call $\Phi:M\times\mathbb{R}\rightarrow M$ the flow of $V$ and $L_{p}$ the
orthogonal leaf of $V$ through $p\in M$. Since $\Phi_{t}$ is an isometry which
preserves $V$, it holds $\Phi_{t}(L_{p})=L_{\Phi_{t}(p)}$ for all
$t\in\mathbb{R}$, i.e., it is a foliated map for each $t\in\mathbb{R}$. Therefore it is easy to show
that the restriction $\Phi:L\times\mathbb{R}\rightarrow M$, where $L$ is a
fixed leaf, is a local diffeomorphism and so $\Phi\left(  L\times
\mathbb{R}\right)  =\cup_{t\in\mathbb{R}}\Phi_{t}(L)$ is an open subset of $M$.

If $p\not \in \cup_{t\in\mathbb{R}}\Phi_{t}(L)$, we can show in the same way
as before that $\cup_{t\in\mathbb{R}}\Phi_{t}(L_{p})$ is an open neighborhood
of $p$, which is contained in the complementary of $\cup_{t\in\mathbb{R}}%
\Phi_{t}(L)$. Since $M$ is supposed connected, $\Phi$ is onto.

Now, take the pull-back metric $\Phi^{\ast}(g)$ which makes $\Phi$ a local
isometry. If $v\in T_{p}L\prec T_{p}M$ then
\[
|(v_{p},0_{t})|=|(\Phi_{t})_{\ast p}(v_{p})|=|v_{p}|=|(\Phi_{0})_{\ast
p}(v_{p})|=|(v_{p},0_{0})|,
\]
where we have taken into account that $\Phi|_{L\times\{0\}}=id$. On the other
hand, if we call $f(p,t)=|V_{\Phi(p,t)}|$ then $f(p,t)=f(p,0)$ and therefore
we can conclude that $\Phi^{\ast}(g)=g|_{L}+\varepsilon f^{2}dt^{2}$, where
$\varepsilon$ is the sign of $V$.

Now we show that it is a covering map. Let $\sigma:[0,1]\rightarrow M$ be a
geodesic and $(x_{0},t_{0})\in L\times\mathbb{R}$ a point such that
$\Phi(x_{0},t_{0})=\sigma(0)$. We must show that there exists a lift
$\alpha:[0,1]\rightarrow L\times\mathbb{R}$ of $\sigma$ through $\Phi$ starting
at $(x_{0},t_{0})$, \cite{ONeill}. There is a geodesic $\alpha:[0,s_{0}%
)\rightarrow L\times\mathbb{R}$, $\alpha(s)=(x(s),t(s))$, such that $\Phi
\circ\alpha=\sigma$ and $\alpha(0)=(x_{0},t_{0})$ because $\Phi$ is a local
isometry. If we suppose $s_{0}<1$, there is a geodesic $(x_{1}(s),t_{1}(s))$
such that $\Phi(x_{1}(s),t_{1}(s))=\sigma(s)$ with $s\in(s_{0}-\delta
,s_{0}+\delta)$. Then in the open interval $(s_{0}-\delta,s_{0})$ it holds
$\Phi(x(s),t(s))=\Phi(x_{1}(s),t_{1}(s))$. Differentiating we get for
$s\in(s_{0}-\delta,s_{0})$ that
\[
\Phi_{\ast_{\left(  x(s),t(s)\right)  }}\left(  t^{\prime}(s)\partial
_{t}+x^{\prime}(s)\right)  =\Phi_{\ast_{\left(  x_{1}(s),t_{1}(s)\right)  }%
}\left(  t_{1}^{\prime}(s)\partial_{t}+x_{1}^{\prime}(s)\right)  ,
\]
and since $\Phi_{*_{(x,t)}}(\partial_t)=V_{\Phi_t(x)}$
\[
\left(  t^{\prime}(s)-t_{1}^{\prime}(s)\right)  V_{\sigma(s)}=\Phi
_{\ast_{\left(  x(s),t(s)\right)  }}\left(  x^{\prime}(s)\right)  -\Phi
_{\ast_{\left(  x_{1}(s),t_{1}(s)\right)  }}\left(  x_{1}^{\prime}(s)\right)
.
\]

But $\Phi_{*_{\left(  x(s),t(s)\right)  }}\left(  x^{\prime}(s)\right)
-\Phi_{*_{\left(  x_{1}(s),t_{1}(s)\right)  }}\left(  x_{1}^{\prime
}(s)\right)  $ is orthogonal to $V$ since $\Phi_{t}$ is a foliated map for each $t\in\mathbb{R}$. Therefore
$t_{1}(s)-t(s)=c\in\mathbb{R}$ and so it exists $\lim_{s\rightarrow s_{0}}
t(s)$ and since $x(s)=\Phi_{-t(s)}(\sigma(s))$ it also exists $\lim
_{s\rightarrow s_{0}} x(s)$ and $\lim_{s\rightarrow s_{0}} \alpha(s)$. Thus
the geodesic $\alpha$ is extendible.

It remains to show that the group of deck transformations acts transitively on
the fibre. Fix $p\in L$ and take $(x_{0},t_{0})\in L\times\mathbb{R}$ such
that $\Phi(x_{0},t_{0})=\Phi(p,0)=p$. Since $\Phi_{-t_{0}}$ is a foliated map, in the sense that it preserves the
foliation given by $V^\perp$, it
follows that $\Phi_{-t_{0}}(L)=L$ and the map $L\times\mathbb{R}\rightarrow
L\times\mathbb{R}$ given by $(x,t)\rightarrow(\Phi_{-t_{0}}(x),t+t_{0})$ is a
deck transformation and takes $(p,0)$ to $(x_{0},t_{0})$.
\end{proof}

\begin{remark}
\label{remarkdglobal} \rm{In the above proposition, we can ensure the
injectivity of $\Phi$, and thus the global decomposition of $M$, supposing
that integral curves of $V$ only intersect each orthogonal leaf one time.
Anyway, although we do not suppose completeness, we can still obtain a local
decomposition.}
\end{remark}

\section{Killing vector fields in static manifolds}

In order to tackle the uniqueness problem, we start studying Killing vector
fields in a static manifold.

\begin{proposition}
\label{Killing} Let $M=L\times_{\varepsilon f}\mathbb{R}$ be a static manifold
and $V=a\partial_{t}+W$, where $a\in C^{\infty}(M)$ and $W\in\mathfrak{X}(M)$
with $W\perp\partial_{t}$. Then $V$ is a Killing vector field if and only if
$V^{t}$ is a Killing vector field on $L$ for each $t\in\mathbb{R}$ and the
following equations hold
\begin{align}
\frac{d}{dt}V^{t}  &  =-\varepsilon f^{2}\nabla^{L}a^{t},\label{eqKilling1}\\
V^{t}(\ln f)  &  =-a_{t}. \label{eqKilling2}%
\end{align}

\end{proposition}

\begin{proof}
Since $L\times\{t\}$ is a geodesic hypersurface of $M$, it is straightforward
that $V^{t}$ is a Killing field on $L$ for all $t\in\mathbb{R}$. We state
equations (\ref{eqKilling1}) and (\ref{eqKilling2}).

Given $X\in\mathfrak{X}(L)$, from $g(\nabla_{\partial_{t}}V,\tilde
{X})=-g(\nabla_{\tilde{X}}V,\partial_{t})$, we get
\begin{align*}
&  ag(\nabla_{\partial_{t}}\partial_{t},\tilde{X})+g(\nabla_{\partial_{t}%
}W,\tilde{X})=\\
&  -\varepsilon f^{2}\tilde{X}(a)-ag(\nabla_{\tilde{X}}\partial_{t}%
,\partial_{t})-g(\nabla_{\tilde{X}}W,\partial_{t}).
\end{align*}
Since $\partial_{t}$ is Killing and $g(\nabla_{\tilde{X}}W,\partial
_{t})=-g(W,\nabla_{\tilde{X}}\partial_{t})=0$, the above simplifies to
\[
g(\nabla_{\partial_{t}}W,\tilde{X})=-\varepsilon f^{2}\tilde{X}(a).
\]
We compute each term at a point $(p_{0},t_{0})\in L\times\mathbb{R}$. The
projection $\pi:L\times\{t_{0}\}\rightarrow L$ is an isometry, so we have
\begin{align*}
&  g(\nabla_{\partial_{t}}W,\tilde{X})_{(p_{0},t_{0})}=\partial_{t}%
|_{(p_{0},t_{0})}g(W,\tilde{X})=\frac{d}{dt}g(W,\tilde{X})_{(p_{0}%
,t)}|_{t=t_{0}}=\\
&  \frac{d}{dt}g_{L}\left(  \pi_{\ast(p_{0},t)}(W),X_{p_{0}}\right)
|_{t=t_{0}}=\frac{d}{dt}g_{L}\left(  W_{p_{0}}^{t},X_{p_{0}}\right)
|_{t=t_{0}}=\\
&  g_{L}\left(  \frac{d}{dt}W_{p_{0}}^{t}|_{t=t_{0}},X_{p_{0}}\right)
=g_{L}\left(  \frac{d}{dt}V_{p_{0}}^{t}|_{t=t_{0}},X_{p_{0}}\right)  .
\end{align*}
On the other hand, if $\gamma$ is an integral curve of $X$, then
\begin{align*}
\tilde{X}(a)_{(p_{0},t_{0})}  &  =\frac{d}{dt}a(\gamma(t),t_{0})|_{t=t_{0}%
}=\frac{d}{dt}a^{t_{0}}(\gamma(t))|_{t=t_{0}}\\
&  =X_{p_{0}}(a^{t_{0}})=g(\nabla a^{t_{0}},X_{p_{0}}).
\end{align*}
Therefore, $\frac{d}{dt}V^{t}=-\varepsilon f^{2}\nabla^{L}a^{t}$ for all
$t\in\mathbb{R}$.

From $g(\nabla_{\partial_{t}}V,\partial_{t})=0$ it follows $g(\nabla
_{\partial_{t}}W,\partial_{t})=-\varepsilon f^{2}a_{t}$. But
\begin{align*}
g(\nabla_{\partial_{t}}W,\partial_{t})  &  =-g(W,\nabla_{\partial_{t}}%
\partial_{t})=\varepsilon fg(W,\widetilde{\nabla f})\\
&  =\varepsilon fW^{t}(f)=\varepsilon fV^{t}(f),
\end{align*}
and therefore $V^{t}(\ln f)=-a_{t}$ for all $t\in\mathbb{R}$. The ``only if''
part is clear.
\end{proof}

Observe that $\frac{d}{dt}V^{t}$ is the vectorial derivative of $V^{t}$, thus it
is also a Killing vector field.

\begin{theorem}
\label{teoKilling} Let $M=L\times_{\varepsilon f}\mathbb{R}$ be a static
manifold with $L$ complete. If $V\in\mathfrak{X}(M)$ is a Killing vector
field, then one of the following holds.

\begin{enumerate}
\item $V=(a_{1}t+a_{2})\partial_{t}+\widetilde{W}$ where $a_{1},a_{2}%
\in\mathbb{R}$ and $W\in\mathfrak{X}(L)$ is a Killing vector field with $W(\ln
f)=-a_{1}$.

\item $L$ decomposes as a static manifold $\left(  N\times\mathbb{R}%
,g_{N}+\lambda(x)^{2}ds^{2}\right)  $ and $f(x,s)=\lambda(x)c(s)$ for certain
$c\in C^{\infty}(\mathbb{R})$.
\end{enumerate}
\end{theorem}

\begin{proof}
Decompose $V=a\partial_{t}+W$ where $a\in C^{\infty}(M)$ and $W\perp
\partial_{t}$. Using Proposition \ref{Killing}, $V^{t}$ is a Killing field on
$L$ and $\frac{d}{dt}V^{t}=-\varepsilon f^{2}\nabla^{L}a^{t}$. Therefore,
$\frac{d}{dt}V^{t}$ is a Killing and irrotational vector field for each
$t\in\mathbb{R}$ .

Fix $t\in\mathbb{R}$, take $\gamma$ a geodesic in $L$ with $\gamma(0)=p$ and
call $y(s)=a^{t}(\gamma(s))$. Then $y_{s}=g_{L}(\nabla^{L}a^{t},\gamma
^{\prime})$ and
\begin{align*}
y_{ss}  &  =g_{L}(\nabla_{\gamma^{\prime}}^{L}\nabla^{L}a^{t},\gamma^{\prime
})=\varepsilon\frac{2(f\circ\gamma)_{s}}{(f\circ\gamma)^{3}}g_{L}\left(
\frac{d}{dt}V^{t},\gamma^{\prime}\right)  =-\frac{2(f\circ\gamma)_{s}}%
{f\circ\gamma}g_{L}(\nabla^{L}a^{t},\gamma^{\prime})\\
&  =-2(\ln f\circ\gamma)_{s}y_{s}.
\end{align*}

Therefore, $y_{s}(s)=y_{s}(0)\left(  \frac{f(p)}{f(\gamma(s))}\right)  ^{2}$.
If $\frac{d}{dt}V_{p}^{t}=0$ then $\nabla^L a^t=0$ at $p$, hence $y_s(0)=0$ and
thus $y_{s}\equiv0$. Therefore $a^{t}$ is constant in a
neighborhood of $p$, which implies that $\frac{d}{dt}V^{t}=0$ in this same
neighborhood. Since $L$ is supposed connected, for each $t\in\mathbb{R}$ there
are two possibilities: $\frac{d}{dt}V^{t}\equiv0$ or it does not have any zero.

If $\frac{d}{dt}V^{t}\equiv0$ for all $t\in\mathbb{R}$, then $V=a\partial
_{t}+\widetilde{W}$ where $W\in\mathfrak{X}(L)$ is a Killing field. Moreover,
by equation (\ref{eqKilling1}), $a$ only depends on $t$ and equation
(\ref{eqKilling2}) implies that there is a constant $a_{1}\in\mathbb{R}$ with
$a_{t}=-W(\ln f)=a_{1}$, from which the first assertion follows.

Suppose on the contrary that there is some $t_{0}\in\mathbb{R}$ such that
$\frac{d}{dt}V^{t_{0}}$ has not zeros. If $\alpha:\mathbb{R}\rightarrow L$ is
an integral curve of $\frac{d}{dt}V^{t_{0}}$, then $\frac{d}{ds}a^{t_{0}%
}(\alpha(s))\neq0$ for all $s\in\mathbb{R}$. Therefore, since $a$ is constant
through the orthogonal leaves of $\frac{d}{dt}V^{t_{0}}$, $\alpha$ only
intersects them one time, which implies that $L$ decomposes as a static
manifold $N\times_{\lambda}\mathbb{R}$, where $\partial_{s}=\frac{d}%
{dt}V^{t_{0}}$ (see Remark \ref{remarkdglobal}). Since $\partial
_{s}=-\varepsilon f^{2}\nabla a^{t_{0}}$, then $a^{t_{0}}$ only depends on $s$
and a direct computation gives us $\nabla a^{t_{0}}=\frac{a_{s}^{t_{0}}%
}{\lambda^{2}}\partial_{s}$. Replacing in the above equation, $f^{2}%
=-\varepsilon\frac{\lambda^{2}}{a_{s}^{t_{0}}}$ and thus $f(x,s)=\lambda
(x)c(s)$, where $c(s)=\sqrt{\frac{-\varepsilon}{a_{s}^{t_{0}}}}$.
\end{proof}

\begin{remark}
\rm{Recall that in the above theorem, as in other results in this paper, when $L$
is one dimensional the factor $N$ is a point and therefore can be removed.}
\end{remark}

\begin{corollary}
\label{basecompacta} If $M=L\times_{\varepsilon f}\mathbb{R}$ is a static
manifold with $L$ compact, then any Killing vector field is of the form
$a\partial_{t}+\widetilde{W}$ where $a\in\mathbb{R}$ and $W\in\mathfrak{X}(L)$ is
a Killing vector field with $W(f)=0$.
\end{corollary}

\begin{proof}
Since $L$ is compact, only the first case of the above theorem holds.
Moreover, $f$ must have a critical point and so $a_{1}=0$.
\end{proof}

A condition on the lightlike sectional curvature also gives us information
about Killing vector fields.

\begin{corollary}
\label{coroKilling} Let $M=L\times_{-f}\mathbb{R}$ be a static space with $L$
complete and dimension greater than one. If there is a point $(p_{0},t_{0})\in
M$ such that $\mathcal{K}(\Pi)\neq0$ for any degenerate plane $\Pi$ of
$T_{(p_{0},t_{0})}M$, then any Killing vector field is of the form
$(a_{1}t+a_{2})\partial_{t}+\widetilde{W}$, where $a_{1},a_{2}\in\mathbb{R}$
and $W\in\mathfrak{X}(L)$ is a Killing vector field with $W(\ln f)=-a_{1}$.
\end{corollary}

\begin{proof}
Suppose that $L$ can be decomposed as $\left(  N\times\mathbb{R},g_{N}%
+\lambda(x)^{2}ds^{2}\right)  $, being $p_{0}$ identified with $(x_{0},s_{0}%
)$, and $f(x,s)=\lambda(x)c(s)$. Take a unitary vector $v\in T_{x_{0}}N$ and
$\Pi$ the degenerate plane spanned by $v$ and $u=\frac{1}{\lambda(x_{0}%
)}\partial_{s}-\frac{1}{f(s_{0},x_{0})}\partial_{t}$. We know that
$\mathcal{K}_{u}(\Pi)=K^{L}\left(  span(v,\partial_{s})\right)  +\frac
{g(\nabla_{v}\widetilde{\nabla f},v)}{f}$, but being $L$ also a static
manifold,
\begin{align*}
K^{L}\left(  span(v,\partial_{s})\right)   &  =-\frac{g(\nabla_{v}%
\widetilde{\nabla\lambda},v)}{\lambda},\\
g(\nabla_{v}\widetilde{\nabla f},v)  &  =cg(\nabla_{v}\widetilde{\nabla
\lambda},v).
\end{align*}
Thus $\mathcal{K}_{u}(\Pi)=0$, which is a contradiction. Applying Theorem
\ref{teoKilling} we get the conclusion.
\end{proof}

\section{Standard static vector fields}

In this section, we show that if a static manifold $M$ admits two different
static decompositions, then the base itself is a static manifold and $M$ can
be viewed as a special type of warped product. A standard static vector field
in a Lorentzian manifold will be always supposed timelike.

\begin{proposition}
\label{static} Let $M=L\times_{\varepsilon f}\mathbb{R}$ be a static manifold.
If $V=a\partial_{t}+W$, where $a\in C^{\infty}(M)$ and $W\perp\partial_{t}$,
is a static vector field, then $V^{t}$ is a static vector field on $L$ for
each $t\in\mathbb{R}$. Moreover, if $\dim L\geq2$, in the open set $\{p\in
L:a^{t}\neq0,V^{t}\neq0\}$ it holds
\begin{equation}
X\left(  \ln(a^{t}f)\right)  =X\left(  \ln\sqrt{g(V^{t},V^{t})}\right)
\label{eqirrotacional}%
\end{equation}
for all $X\in\mathfrak{X}(L)$ with $X\perp V^{t}$.
\end{proposition}

\begin{proof}
The first assertion follows easily because $L\times\{t\}$ are geodesic
hypersurfaces. Suppose that $dimL\geq2$ and take $\xi=-\varepsilon
g(W,W)\partial_{t}+af^{2}W$, which is a vector field orthogonal to $V$. Since
$V$ is static, $g(\nabla_{\tilde{X}}V,\xi)=0$ for all $X\in\mathfrak{X}(L)$
with $X\perp W$. But
\begin{align*}
g(\nabla_{\tilde{X}}V,\xi)  &  =\left(  \tilde{X}(a)+a\tilde{X}(\ln f)\right)
g(\partial_{t},\xi)+g(\nabla_{\tilde{X}}W,\xi)=\\
&  =-\left(  \tilde{X}(a)+a\tilde{X}(\ln f)\right)  f^{2}g(W,W)+\frac{af^{2}%
}{2}\tilde{X}\left(  g(W,W)\right) \\
&  =-\left(  X(a^{t})+a^{t}X(\ln f)\right)  g(V^{t},V^{t})f^{2}+\frac
{a^{t}f^{2}}{2}X\left(  g(V^{t},V^{t})\right)  .
\end{align*}

Therefore $\left(  X(a^{t})+a^{t}X(\ln f)\right)  g(V^{t},V^{t})=\frac{a^{t}%
}{2}X(g(V^{t},V^{t}))$. Where $a^{t}\neq0$ and $V^{t}\neq0$, we can write
\[
X\left(  \ln(a^{t}f)\right)  =X\left(  \ln\sqrt{g(V^{t},V^{t})}\right)  .
\]

\end{proof}

\begin{proposition}
\label{descompL} Let $M=L\times_{\varepsilon f}\mathbb{R}$ be a static
manifold with $L$ complete and $V$ a static vector field linearly independent
to $\partial_{t}$ at some point. Then it holds the following.

\begin{enumerate}
\item For each $t\in\mathbb{R}$, $V^{t}$ is identically zero or it does not
have zeros. In fact, there exists a dense open subset $\Theta\subset
\mathbb{R}$ such that the second statement holds for all $t\in\Theta$.

\item If moreover $V$ is standard, then $V^{t}$ is standard static in $L$ for each
$t\in\Theta$. So, fixed $t\in\Theta$, $L$ decomposes as $\left(
N\times\mathbb{R},g_{N}+\lambda(x)^{2}ds^{2}\right)  $ where $V^{t}$ is
identified with $\partial_{s}$.
\end{enumerate}
\end{proposition}

\begin{proof}
\

\begin{enumerate}
\item Fixed $t\in\mathbb{R}$, call $A=\{p\in L:\text{ $V$ and $\partial_{t}$
are l.i. at $(p,t)$ }\}$ and $B=A^{c}$. It is clear that $A$ is open and since
the orthogonal leaf of $V$ and $\partial_{t}$ are geodesic, $B$ is also open.
Therefore, $A=L$ and thus $V^{t}$ does not have zeros or $B=L$ and
$V^{t}\equiv0 $. Now, call $\Theta=\{t\in\mathbb{R}:V^{t}\text{ does not have
zeros}\}$, which obviously is open. If $V^{t}\equiv0$ for all $t\in
(-\delta,\delta)$, then $V$ and $\partial_{t}$ are linearly dependent in
$L\times(-\delta,\delta) $, but since they are Killing vector fields, they
must be linearly dependent in the whole $M$, which is a contradiction.
Therefore $\Theta$ is dense.

\item Given $t\in\Theta$, using the above point and Proposition \ref{static},
we know that $V^{t}$ is a static vector field without zeros in $L$. We show
that it gives rise to a global decomposition. Call $F_{(p,t)}$ the orthogonal
leaf of $V$ through $(p,t)$ and $N$ the orthogonal leaf of $V^{t}$ through $p$,
which is inside $\pi\left(  L\times\{t\}\cap F_{(p,t)}\right)  $. Take
$\alpha$ an integral curve of $V^{t}$ with $\alpha(0)\in N$ and suppose that
there is $s>0$ with $\alpha(s)\in N$. Since $V$ is standard, there is a global
projection $P:M\rightarrow\mathbb{R}$ such that it is constant through the
orthogonal leaves of $V$ and $P_{\ast_{(p,t)}}(v)$ gives the component in the
direction of $V_{(p,t)}$ of any vector $v\in T_{(p,t)}M$. If we call
$\gamma(s)=(\alpha(s),t)$, then $P(\gamma(s))$ has a critical point $s_{1}%
\in(0,s)$, because $\gamma(0),\gamma(s)\in F_{(p,t)}$. But then
\begin{align*}
&  g\left(  V_{\alpha(s_{1})}^{t},V_{\alpha(s_{1})}^{t}\right)  =g\left(
\gamma^{\prime}(s_{1}),V_{\gamma(s_{1})}\right)  =\\
&  P_{\ast_{\gamma(s_{1})}}\left(  \gamma^{\prime}(s_{1})\right)  g\left(
V_{\gamma(s_{1})},V_{\gamma(s_{1})}\right)  =0,
\end{align*}
which is a contradiction. Using Remark \ref{remarkdglobal}, $L$ can be
decomposed as
\[
\left(  N\times\mathbb{R},g_{N}+\lambda(x)^{2}ds^{2}\right)  ,
\]
where $\partial_{s}$ is identified with $V^{t}$.
\end{enumerate}
\end{proof}

As a consequence of the above proposition, we can prove the main result of
\cite{Sanchez1} in a different way.

\begin{theorem}
\label{MiguelSenovilla}Let $M=L\times_{\varepsilon f}\mathbb{R}$ be a static
manifold with $L$ compact. Then any other standard static vector field is
proportional to $\partial_{t}$ and so $M$ admits a unique decomposition as a
static manifold.
\end{theorem}

The following Proposition will be the key to prove our main results.

\begin{proposition}
\label{clasifica}Let $M=L\times_{\varepsilon f}\mathbb{R}$ be a static
manifold with $L$ complete and $V$ a standard static vector field linearly
independent to $\partial_{t}$ at a point $(p_{0},t_{0})$. If $V^{t}$ is
proportional to $V^{t_{0}}$ for all $t\in\mathbb{R}$, then $M$ is isometric to
one of the following.

\begin{enumerate}
\item $(N\times\mathbb{R}^{2},g_{N}+\lambda(x)^{2}ds^{2}+\varepsilon
f(x)^{2}dt^{2})$ and $V=\partial_{s}$.

\item $\left(  N\times\mathbb{H}_{[\varepsilon]}^{2}(r),g_{N}+\lambda
(x)^{2}(ds^{2}+\varepsilon\cosh^{2}(rs)dt^{2})\right)  $ and
\[
V=\left(  -\frac{\varepsilon}{r}h_{t}(t)\tanh(rs)+\gamma\right)  \partial
_{t}+h(t)\partial_{s},
\]
where $h(t)=\alpha\sin(rt+\beta)$ if $\varepsilon=-1$ or $h(t)=\alpha
e^{rt}+\beta e^{-rt}$ if $\varepsilon=1$ and $\alpha,\beta,\gamma\in
\mathbb{R}$.

\item $\left(  N\times\widehat{\mathbb{H}}_{[\varepsilon]}^{2}(r),g_{N}%
+\lambda(x)^{2}(ds^{2}+\varepsilon e^{2rs}dt^{2})\right)  $ and
\[
V=\left(  \frac{\varepsilon\alpha}{2r}e^{-2rs}-\frac{r\alpha}{2}t^{2}-r\beta
t+\gamma\right)  \partial_{t}+(\alpha t+\beta)\partial_{s},
\]
where $\alpha,\beta,\gamma\in\mathbb{R}$.

\item $\left(  N\times\mathbb{R}_{[\varepsilon]}^{2},g_{N}+\lambda
(x)^{2}(ds^{2}+\varepsilon dt^{2})\right)  $ and $V=\gamma\partial
_{t}+\partial_{s}$, where $\gamma\in\mathbb{R}$.
\end{enumerate}
\end{proposition}

\begin{proof}
Suppose that $V^{t}=h(t) V^{t_{0}}$ for some $h\in C^{\infty}(\mathbb{R})$
with $h(t_{0})=1$. Proposition \ref{descompL} ensures that $L$ decomposes as
$\left(  N\times\mathbb{R},g_{N}+\lambda(x)^{2}ds^{2}\right)  $ where
$\partial_{s}$ is identified with $V^{t_{0}}$.

If $X\in\mathfrak{X}(N)$, using equation (\ref{eqKilling1}) of Proposition
\ref{Killing} we get $g(X,\nabla a^{t})=0$ and thus $a$ only depends on $s$
and $t$. On the other hand, multiplying by $V^{t_{0}}$ in equation
(\ref{eqKilling1}) we get
\begin{equation}
a_{s}(t,s)=-\varepsilon\frac{h_{t}(t)\lambda(x)^{2}}{f(x,s)^{2}},
\label{eq1clasifica}%
\end{equation}

and equation (\ref{eqKilling2}) of Proposition \ref{Killing} can be written
as
\begin{equation}
a_{t}(s,t)=-h(t)(\ln f)_{s}. \label{eq2clasifica}%
\end{equation}

Now we consider two possibilities:

1) $a\equiv0$. Then above equations give us that $h\equiv1$, and $f$ only
depends on $x$, i.e., $M$ is isometric to $(N\times\mathbb{R}^{2}%
,g_{N}+\lambda(x)^{2}ds^{2}+\varepsilon f(x)^{2}dt^{2})$ and $V=\partial_{s}$.

2) There is a point $(s_{0},t_{0})$ with $a(s_{0},t_{0})\neq0$. Take
\[
A=\{(s,t)\in\mathbb{R}^{2}:a(s,t)\neq0\}.
\]
Equation (\ref{eqirrotacional}) of Proposition \ref{static} reduces to $X(\ln
f)=X(\ln\lambda)$ in $N\times A$ for all $X\in\mathfrak{X}(N)$, which implies
that there is certain function $c$ such that $f(x,s)=\lambda(x)c(s)$ for all
$(x,s,t)\in N\times A$.

Take $B=\left(  \overline{A}\right)  ^{c}$. If $B=\emptyset$, then by
continuity $f(x,s)=\lambda(x)c(s)$ for all $(x,s,t)\in N\times\mathbb{R}^{2}$.
If $B\not =\emptyset$, then $a\equiv0$ in $N\times B$ and so $f$ only depends
on $x$, i.e. $f(x,s)=F(x)$ for all $(x,s,t) \in N\times B$ where $F$ is
certain function. Since $\lambda(x) c(s)=F(x)$ for all $(x,s,t)\in Fr(N\times
A)=N\times Fr(A)$, it it easy to show that $c$ can be extended to the whole
$\mathbb{R}$ and, with this extension, it holds $f(x,s)=\lambda(x)c(s)$ for
all $(x,s)\in N\times\mathbb{R}^{2}$.

If we call $S=\mathbb{R}\times_{\varepsilon c}\mathbb{R}$, then $M$ is the
warped product $N\times_{\lambda}S$. We now classify this surface $S$.
Equations (\ref{eq1clasifica}) and (\ref{eq2clasifica}) reduce to
\begin{align}
a_{s}  &  =-\varepsilon\frac{h_{t}(t)}{c(s)^{2}},\nonumber\\
a_{t}  &  =-h(t)\left(  \ln c(s)\right)  _{s},\nonumber
\end{align}
and using the Schwarz's Theorem we get the differential equations
\begin{align}
c_{ss}(s)c(s)-c_{s}(s)^{2}  &  =k,\label{eq3clasifica}\\
h_{tt}(t)  &  =\varepsilon kh(t)\nonumber
\end{align}
for some constant $k\in\mathbb{R}$. The solutions of (\ref{eq3clasifica}) are

\begin{itemize}
\item $c(s)=\frac{\sqrt{-k}}{r}\sinh(rs+b)$ or $c(s)=\frac{\sqrt{-k}}{r}%
\sin(rs+b)$ if $k<0$.

\item $c(s)=e^{rs+b}$ if $k=0$.

\item $c(s)=\frac{\sqrt{k}}{r}\cosh(rs+b)$ if $k>0$.
\end{itemize}

Since $c(s)>0$ for all $s\in\mathbb{R}$ we should discard the case $k<0$. In
the case $k>0$, solving the above differential equations, there are
$\alpha,\beta,\gamma\in\mathbb{R}$ such that
\[
V=\left(  \frac{r\alpha}{\sqrt{k}}\cos(\sqrt{k}t+\beta)\tanh(rs+b)+\gamma
\right)  \partial_{t}+\alpha\sin(\sqrt{k}t+\beta)\partial_{s}%
\]
if $\varepsilon=-1$ or
\[
V=\left(  -\frac{r}{\sqrt{k}}\left(  \alpha e^{\sqrt{k}t}-\beta e^{-\sqrt{k}%
t}\right)  \tanh(rs+b)+\gamma\right)  \partial_{t}+\left(  \alpha e^{\sqrt
{k}t}+\beta e^{-\sqrt{k}t}\right)  \partial_{s}%
\]
if $\varepsilon=1$. Now, we obtain point (2) rescaling with $\Phi
:\mathbb{R}^{2}\rightarrow\mathbb{R}^{2}$ given by $\Phi(s,t)=(s+\frac{b}%
{r},\frac{\sqrt{k}}{r}t)$.

If $k=0$, we can suppose $b=0$ rescaling $s$, and thus $c(s)=e^{rs}$ with
$r\neq0$ or $c(s)=1$. In the first case, $S=\widehat{\mathbb{H}}%
_{[\varepsilon]}^{2}(r)$ and
\[
V=\left(  \frac{\varepsilon\alpha}{2r}e^{-2rs}-\frac{r\alpha}{2}t^{2}-r\beta
t+\gamma\right)  \partial_{t}+(\alpha t+\beta)\partial_{s}.
\]

In the second case $S=\mathbb{R}_{[\varepsilon]}^{2}$, $a(s,t)=-\varepsilon
\alpha s+\gamma$ and $h(t)=\alpha t+\beta$. But $\alpha=0$ and $\beta=1$
because $V$ does not have zeros and $h(t_{0})=1$.
\end{proof}

\begin{remark}
\rm{The case $N\times_{\lambda}\mathbb{H}^{2}(r)$ and $N\times_{\lambda}%
\widehat{\mathbb{H}}^{2}(r)$ are equivalent since $\mathbb{H}^{2}(r)$ and
$\widehat{\mathbb{H}}^{2}(r)$ are isometric spaces. However, $N\times
_{\lambda}\mathbb{H}_{[-1]}^{2}(r)$ and $N\times_{\lambda}\widehat{\mathbb{H}%
}_{[-1]}^{2}(r)$ are not equivalent because $\mathbb{H}_{[-1]}^{2}(r)$ is
complete and $\widehat{\mathbb{H}}_{[-1]}^{2}(r)$ is not.}
\end{remark}

As an immediate consequence we obtain the following.

\begin{corollary}
\label{clasificadim2} Let $M=L\times_{\varepsilon f}\mathbb{R}$ be a complete
two dimensional static manifold. If there exists $V$ a non-identically zero
Killing vector field linearly independent to $\partial_{t}$ at some point,
then $M$ is isometric to $\mathbb{R}_{[\varepsilon]}^{2}$ or $\mathbb{H}%
_{[\varepsilon]}^{2}(r)$.
\end{corollary}

\begin{proof}
Since $dimL=1$, $V$ is also irrotational, $V^{t}$ is linearly dependent to
a fixed $V^{t_{0}}$ for all $t\in\mathbb{R}$ and the proof of Proposition
\ref{clasifica} works with $N$ reduced to a point, although $V$ is not
necessarily standard and, maybe, with zeros.
\end{proof}

Now, we show that a manifold with more than one static decomposition is a
particular type of warped product.

\begin{proposition}
\label{teoclasifica} Let $M=L\times_{\varepsilon f}\mathbb{R}$ be a static
manifold with $L$ complete. If there exists $V$ a standard static vector field
linearly independent to $\partial_{t}$ at some point, then $M$ decomposes as
$(N\times\mathbb{R}^{2},g_{N}+\lambda(x)^{2}ds^{2}+\varepsilon f(x)^{2}%
dt^{2})$ and $V=\partial_{s}$ or $L$ decomposes as $\left(  N\times
\mathbb{R},g_{N}+\lambda(x)^{2}ds^{2}\right)  $ and $f(x,s)=\lambda(x)c(s)$,
i.e., $M$ is the warped product $N\times_{\lambda}\left(  \mathbb{R}\times
_{c}\mathbb{R}\right)  $.
\end{proposition}

\begin{proof}
Using Theorem \ref{teoKilling}, $M$ decomposes as a warped product
$N\times_{\lambda}\left(  \mathbb{R}\times_{c}\mathbb{R}\right)  $ or
$V=(a_{1}t+a_{2})\partial_{t}+\widetilde{W}$ where $a_{1},a_{2}\in\mathbb{R}$
and $W\in\mathfrak{X}(L)$. But in this last case, we can apply Proposition
\ref{clasifica} to obtain that $M$ may also be decomposed as $(N\times
\mathbb{R}^{2},g_{N}+\lambda(x)^{2}ds^{2}+\varepsilon f(x)^{2}dt^{2})$ and
$V=\partial_{s}$ or again as $N\times_{\lambda}\left(  \mathbb{R}\times
_{c}\mathbb{R}\right)  $.
\end{proof}

\begin{remark}
\rm{As we said in the introduction, euclidean, Minkowski, hyperbolic and a portion of the anti De Sitter spaces have different static decompositions and therefore they must
fulfill the above proposition. But this can be easily checked taking into
account that the hyperbolic space can be viewed as
\[
\mathbb{H}^{n}=\left(  \mathbb{R}^{n},dx_{1}^{2}+e^{2x_{1}}\left(  \sum
_{i=2}^{n}dx_{i}^{2}\right)  \right)
\]
and the above mentioned portion of anti De Sitter as
\[
\left(  \mathbb{R}^{n},dx_{1}^{2}+e^{2x_{1}}\left(
\sum_{i=2}^{n-1}dx_{i}^{2}-dx_{n}^{2}\right)  \right)  .
\]}

\end{remark}

\begin{corollary}
\label{productodirecto} Let $M=L\times_{\varepsilon}\mathbb{R}$ be a static
manifold with $L$ complete and with constant warping function. If $V$ is a standard
static vector field which is linearly independent to $\partial_{t}$ at some
point, then $M$ is isometric to $(N\times\mathbb{R}^{2},g_{N}+\lambda
(x)^{2}ds^{2}+\varepsilon dt^{2})$ where $V=\partial_{s}$ or to a direct
product $N\times\mathbb{R}_{[\varepsilon]}^{2}$.
\end{corollary}

\section{Main results}

We are already able to classify, under a curvature hypothesis, manifolds with
more than one static decomposition. We start assuming that the base has a
point with positive curvature.

\begin{theorem}
\label{ClasificaCurvatura} Let $M=L\times_{\varepsilon f}\mathbb{R}$ be a
static manifold with $L$ complete and dimension greater than one. Suppose that
there exists a standard static vector field $V$ linearly independent to
$\partial_{t}$ at some point and there is $p\in L$ with $K^{L}(\Pi)>0$ for any
plane $\Pi$ of $T_{p}L$. Then $M$ is isometric to

\begin{enumerate}
\item $(N\times\mathbb{R}^{2},g_{N}+\lambda(x)^{2}ds^{2}+\varepsilon
f(x)^{2}dt^{2})$ and $V=\partial_{s}$.

\item $\left(  N\times\mathbb{H}_{[\varepsilon]}^{2}(r),g_{N}+\lambda
(x)^{2}(ds^{2}+\varepsilon\cosh^{2}(rs)dt^{2})\right)  $ and
\[
V=\left(  -\frac{\varepsilon}{r}h_{t}(t)\tanh(rs)+\gamma\right)  \partial
_{t}+h(t)\partial_{s},
\]
where $h(t)=\alpha\sin(rt+\beta)$ if $\varepsilon=-1$ or $h(t)=\alpha
e^{rt}+\beta e^{-rt}$ if $\varepsilon=1$ and $\alpha,\beta,\gamma\in
\mathbb{R}$.

\item $\left(  N\times\widehat{\mathbb{H}}_{[\varepsilon]}^{2}(r),g_{N}%
+\lambda(x)^{2}(ds^{2}+\varepsilon e^{2rs}dt^{2})\right)  $ and
\[
V=\left(  \frac{\varepsilon\alpha}{2r}e^{-2rs}-\frac{r\alpha}{2}t^{2}-r\beta
t+\gamma\right)  \partial_{t}+(\alpha t+\beta)\partial_{s},
\]
where $\alpha,\beta,\gamma\in\mathbb{R}$.

\item $\left(  N\times\mathbb{R}_{[\varepsilon]}^{2},g_{N}+\lambda
(x)^{2}(ds^{2}+\varepsilon dt^{2})\right)  $ and $V=\gamma\partial
_{t}+\partial_{s}$, where $\gamma\in\mathbb{R}$.
\end{enumerate}
\end{theorem}

\begin{proof}
We proceed by induction over $dimL$. Suppose first that $dimL=2$. Using
Proposition \ref{descompL}, $V^{t}$ is standard static for $t$ in a dense open
set $\Theta$. Fix $t_{0}\in\Theta$ and decompose $L$ as $\left(
N\times\mathbb{R},g_{N}+\lambda(x)^{2}ds^{2}\right)  $ where $\partial
_{s}=V^{t_{0}}$. If there is a $t_{1}\in\Theta$ with $V^{t_{1}}$ and
$\partial_{s}$ linearly independent at some point, then Corollary
\ref{clasificadim2} ensures $L=\mathbb{H}^{2}$ or $L=\mathbb{R}^{2}$, which
contradicts the curvature hypothesis. Therefore $V^{t}$ is linearly dependent
to $\partial_{s}$ for all $t\in\mathbb{R}$ and Proposition \ref{clasifica}
proves the statement.

Now, assuming the statement for $dimL=n-1$, we will prove it for $dimL=n$. As
before, $L=\left(  N\times\mathbb{R},g_{N}+\lambda(x)^{2}ds^{2}\right)  $,
where $\partial_{s}=V^{t_{0}}$. If there is $t\in\Theta$ with $V^{t}$ linearly
independent to $\partial_{s}$ at some point, applying the induction
hypothesis, $L$ is isometric to a warped product $S\times_{\mu}\mathbb{H}%
_{[\varepsilon]}^{2}(r)$, $S\times_{\mu}\widehat{\mathbb{H}}_{[\varepsilon
]}^{2}(r)$, $S\times_{\mu}\mathbb{R}_{[\varepsilon]}^{2}$ or to $(S\times
\mathbb{R}^{2},g_{S}+\mu(z)^{2}du^{2}+\lambda(z)^{2}ds^{2})$ where
$V^{t}=\partial_{u}$ and $V^{t_{0}}=\partial_{s}$. In the first three cases,
any tangent plane to the fibre has nonpositive curvature and should be
discarded. In the last case, $V^{t}$ is orthogonal to $V^{t_{0}}$ and we can
obtain a contradiction using the continuity of $V^{t}$ respect to $t$.
Therefore, $V^{t}$ must be linearly dependent to $\partial_{s}$ for all
$t\in\mathbb{R}$ and applying Proposition \ref{clasifica} we get the result.
\end{proof}

We also obtain the same classification if the manifold is Einstein.

\begin{theorem}
\label{Einstein} Let $M=L\times_{\varepsilon f}\mathbb{R}$ be an Einstein
static manifold with $L$ complete and dimension greater than one. If there
exists a standard static vector field $V$ linearly independent to
$\partial_{t} $ at some point, then $M$ is isometric to $\left(
N\times\mathbb{R}^{2},g_{N}+\lambda(x)^{2}ds^{2}+\varepsilon f(x)^{2}
dt^{2}\right)  $ where $V=\partial_{s}$ or to a warped product $N\times
_{\lambda}\mathbb{H}^{2}_{[\varepsilon]}(r)$, $N\times_{\lambda}%
\widehat{\mathbb{H}}^{2}_{[\varepsilon]}(r)$ or $N\times_{\lambda}%
\mathbb{R}^{2}_{[\varepsilon]}$.
\end{theorem}

\begin{proof}
First, note that for any static space $\left(P\times\mathbb{R},g_P+\varepsilon h^2 dt^2\right)$ it holds
\begin{eqnarray}
 Ric(X,Y)&=&Ric^{P}(X,Y)-\frac{1}{h}H_{h}(X,Y) \text{ for $X,Y\in
\mathfrak{X}(P)$},\label{ecRicci1}\\
Ric(\partial_{t},\partial_{t})&=&-\varepsilon h\triangle^{P}h,\label{ecRicci2}
\end{eqnarray}
where $H_h$ is the hessian of $h$ and $\triangle^P$ the laplacian operator in $P$.

Proposition \ref{teoclasifica} says that $M$ is $\left(  N\times\mathbb{R}%
^{2},g_{N}+\lambda(x)^{2}ds^{2}+\varepsilon f(x)^{2}dt^{2}\right)  $ and
$V=\partial_{s}$ or $L$ decomposes as $\left(  N\times\mathbb{R},g_{N}%
+\lambda(x)^{2}ds^{2}\right)$ and the warping function $f$ also decomposes as $f(x,s)=\lambda(x)c(s)$, where $\lambda\in C^\infty(N)$ and $c\in C^\infty (\mathbb{R})$.
Thus, it is enough to analize the second case and show that $c(s)$ is $e^s,\cosh s$ or a positive constant.

Since $M$ is static, Einstein and $dim M\geq 3$, there is $\delta\in\mathbb{R}$ such that equation (\ref{ecRicci1}) can
be written as
\begin{equation}
Ric^{L}(X,Y)=\frac{1}{f}H_{f}(X,Y)+\delta g(X,Y)\text{ for $X,Y\in
\mathfrak{X}(L)$}. \label{RicciL}%
\end{equation}

Since $L$ is also the static manifold $N\times_\lambda\mathbb{R}$, we have from equation (\ref{ecRicci2})
\begin{align*}
Ric^{L}(\partial_{s},\partial_{s})  &  =-\lambda\triangle^{N}\lambda.
\end{align*}

Moreover, since $f(x,s)=\lambda(x)c(s)$,
\begin{align*}
H_{f}(\partial_{s},\partial_{s})  &  =\lambda\left(  cg(\nabla\lambda
,\nabla\lambda)+c_{ss}\right)  .
\end{align*}

Notice that $\partial_s\in\mathfrak{X}(L)$, thus we can replace the above two equations in equation (\ref{RicciL}) to get
\begin{align}
-\lambda\triangle^{N}\lambda &  =g(\nabla\lambda,\nabla\lambda)+\frac{c_{ss}%
}{c}+\delta\lambda^{2}. \label{laplaciano}%
\end{align}

The left hand side and the first and third summands of the right hand side of this equation do not depend on the parameter $s$,
thus we conclude that $\frac{c_{ss}}{c}=a$ for
certain constant $a\in\mathbb{R}$. Moreover, since $L$ is complete and $c$ has not zeros,
necessarily $a=r^{2}\geq0$ and therefore $c(s)=e^{rs+\beta}$ or $c(s)=\alpha
\cosh(rs+\beta)$ with $\alpha\in\mathbb{R}^+,\beta\in\mathbb{R}$. Rescaling we obtain that $M$ is isometric to $N\times
_{\lambda}\mathbb{H}_{[\varepsilon]}^{2}(r)$, $N\times_{\lambda}%
\widehat{\mathbb{H}}_{[\varepsilon]}^{2}(r)$ or $N\times_{\lambda}%
\mathbb{R}_{[\varepsilon]}^{2}$.
\end{proof}

\begin{remark}\rm{If $L$ is the static manifold $N\times_\lambda\mathbb{R}$, then equation $\ref{ecRicci1}$ says
\begin{eqnarray*} 
Ric^{L}(X,Y)  &  =Ric^{N}(X,Y)-\frac{1}{\lambda}H_{\lambda}(X,Y)\text{ for
$X,Y\in\mathfrak{X}(N)$},
\end{eqnarray*}
and since $f(x,s)=\lambda(x)c(s)$,
\begin{eqnarray*}
H_{f}(X,Y)  &  =cH_{\lambda}(X,Y)\text{ for $X,Y\in\mathfrak{X}(N)$}.
\end{eqnarray*}

Replacing the above two equations in (\ref{RicciL}), we get the equation
\begin{eqnarray}
Ric^{N}(X,Y)  &  =\frac{2}{\lambda}H_{\lambda}(X,Y)+\delta g(X,Y)\text{ for
$X,Y\in\mathfrak{X}(N)$}\label{RicciN},
\end{eqnarray}
which will be used in the next theorem.}
\end{remark}

\begin{remark}
\label{remarklocal} \rm{The arguments used in Theorems \ref{ClasificaCurvatura} and
\ref{Einstein} also work locally, so we can avoid the completeness of the base
obtaining similar conclusions. More concretely, we would obtain that a
neighborhood of the point where $\partial_{t}$ and $V$ (being $V$ static but
nonnecessarily standard) are linearly independent, is locally isometric to
$(N\times\mathbb{R}^{2},g_{N}+\lambda(x)^{2}ds^{2}+\varepsilon f(x)^{2}%
dt^{2})$ where $V=\partial_{s}$ or to a warped product $N\times_{\lambda}S$,
where $S$ is a surface of constant curvature.}
\end{remark}

Now, we particularize the above to the case of a (four dimensional) spacetime.
In the following theorem it is shown that a fundamental component of static
Einstein spacetimes with different decompositions are the Riemannian surfaces
\begin{equation}
\frac{1}{k_{1}+\frac{k_{2}}{u}+k_{3}u^{2}}du^{2}+(k_{1}+\frac{k_{2}}{u}%
+k_{3}u^{2})dv^{2}, \label{metrica}%
\end{equation}
where $k_{1},k_{2},k_{3}\in\mathbb{R}$. Observe that surfaces of constant
curvature $1$, $-1$ and $0$ are included in this family for an appropriate
choice of the constants.

\begin{theorem}
\label{Einsteinspacetime} Let $M=L\times_{-f}\mathbb{R}$ be an Einstein static
spacetime. If there exists a timelike static vector field linearly independent
to $\partial_{t}$ at some point, then almost every point in $M$ has a
neighborhood isometric to a direct product of two surfaces with the same
constant curvature or to a warped product $\Gamma\times_{u}S$, where $\Gamma$
is a surface as the one given in (\ref{metrica}) and $S$ has constant curvature.
\end{theorem}

\begin{proof}
As it is said in Remark \ref{remarklocal}, since $L$ is not supposed complete,
from Theorem \ref{Einstein} we only obtain that $M$ is locally a warped
product $N\times_{\lambda}S$, where $S$ has constant curvature (the other case
in the remark should be discarded because a standard static vector field in a
Lorentzian manifold is supposed timelike). Since $\dim N=2$, we have
$Ric^{N}=Kg$, being $K$ the curvature of $N$, and thus equation (\ref{RicciN})
transforms to
\begin{equation}
H_{\lambda}=\frac{\lambda(K-\delta)}{2}g. \label{eqhessiano}%
\end{equation}
Taking trace we get
\begin{equation}
\triangle^{N}\lambda=\lambda(K-\delta) \label{eq1}%
\end{equation}
and replacing in equation (\ref{laplaciano}) we have
\begin{equation}
-\lambda^{2}K=g(\nabla\lambda,\nabla\lambda)+a. \label{eqlambda}%
\end{equation}
If $\lambda$ is constant, then $N$ has constant curvature $K=\frac{-a}%
{\lambda^{2}}$ and rescaling the metric of $S$, $M$ is locally the direct
product of two surfaces with the same constant curvature.

Suppose now that $\lambda$ is not constant. Differentiating the above
equation, $-\lambda dK=(3K-\delta)d\lambda$ and thus $\lambda^{3}%
K=b+\frac{\delta}{3}\lambda^{3}$, for certain $b\in\mathbb{R}$. Now, equation
(\ref{eqlambda}) can be written as
\begin{equation}
-\frac{b}{\lambda}-\frac{\delta}{3}\lambda^{2}=g(\nabla\lambda,\nabla
\lambda)+a. \label{eqlambda2}%
\end{equation}
Equation (\ref{eqhessiano}) means that $\nabla\lambda$ is conformal and, since
it is not constant, its critical points are isolated. Moreover, in a
neighborhood of any point where $\nabla\lambda\neq0$, $N$ is a warped product
$\left(  \mathbb{R}^{2},du^{2}+\lambda_{u}^{2}dv^{2}\right)  $ being
$\partial_{u}=\frac{\nabla\lambda}{|\nabla\lambda|}$, \cite{Montiel, Tashiro,
Kuhnel}. In this coordinates, equation (\ref{eqlambda2}) is
\begin{equation}
\lambda_{u}^{2}=-a-\frac{b}{\lambda}-\frac{\delta}{3}\lambda^{2},
\label{eqlambda3}%
\end{equation}
and reparametrizing with $(u,v)\mapsto(\lambda(u),v)$, $N$ is locally
isometric to
\[
\frac{1}{h(u)}du^{2}+h(u)dv^{2},
\]
where $h(u)=-a-\frac{b}{u}-\frac{\delta}{3}u^{2}$.
\end{proof}

\begin{remark}
\rm{The same conclusion of this Theorem is
obtained in \cite{Tod}, but with slightly different hypothesis. The author
starts with a fixed Riemannian 3-dimensional metric and he supposes that
different Einstein static spacetimes can be constructed using it. Then, he
proves that the metric is locally the one given in (\ref{metrica}). Observe
that if we have a static manifold with two different static decompositions, a
priori, we do not know if the respective bases are isometric.}
\end{remark}

Finally, we show that a condition on the lightlike sectional curvature ensures
the uniqueness of the static decomposition in the Lorentzian case.

\begin{theorem}
\label{Unicidad} Let $M=L\times_{-f}\mathbb{R}$ be a static space with $dim
L\geq2$. If there exists a point $(p,t)\in M$ such that $\mathcal{K}(\Pi
)\neq0$ for any degenerate plane $\Pi$ of $T_{(p,t)}M$, then $M$ admits an
unique decomposition as static space.
\end{theorem}

\begin{proof}
Let $V$ be a (timelike) standard static vector field on $M$ linearly
independent to $\partial_{t}$ at some point and suppose first that $L$ is
complete. Using Proposition \ref{teoclasifica}, $L=\left(  N\times
\mathbb{R},g_{N}+\lambda(x)^{2}ds^{2}\right)  $ and $f(x,s)=\lambda(x)c(s)$
and we can show as in Corollary \ref{coroKilling} that there is a degenerate
plane at $(p,t)$ with zero lighlike sectional curvature, which is a contradiction.

Now, although $L$ is not necessarily complete, Proposition \ref{teoclasifica}
is valid locally and we can still use the above arguments to show that in a
neighborhood $U$ of $(p,t)$ there is a unique static decomposition. Therefore
$V$ and $\partial_{t}$ are linearly dependent in $U$, but being Killing vector
fields, this implies that in fact they are linearly dependent in the whole $M$.
\end{proof}

It would be interesting to have a more accurate classification of static
manifolds admitting more than one standard static vector field at least in
dimension 3 and 4.

\end{document}